\documentclass{article}
\usepackage{amssymb,latexsym,amscd}

\newcommand{\se}[1]{{\section{#1}} {\setcounter{equation}{0}}}
\newtheorem{theorem}{Theorem}[section]
\newtheorem{lm}{Lemma}[section]
\newtheorem{prop}{Proposition}[section]
\newtheorem{de}{Definition}[section]
\newtheorem{co}{Corollary}[section]

\def\k{{K\"{a}hler }}
\def\ke{{K\"{a}hler-Einstein }}
\def\cy{{Calabi-Yau }}

\begin{document}
\hbadness=10000
\title{{\bf H-minimal Lagrangian fibrations in \k manifolds and minimal Lagrangian vanishing tori in \ke manifolds}}
\author{Wei-Dong Ruan\\
Department of Mathematics\\
University of Illinois at Chicago\\
Chicago, IL 60607\\}
\footnotetext{Partially supported by NSF Grant DMS-0104150.}
\maketitle
\begin{abstract}
H-minimal Lagrangian submanifolds in general \k manifolds generalize special Lagrangian submanifolds in Calabi-Yau manifolds. In this paper we will use the deformation theory of H-minimal Lagrangian submanifolds in \k manifolds to construct minimal Lagrangian torus in certain K\"{a}hler-Einstein manifolds with negative first Chern class.
\end{abstract}
%\tableofcontents
\se{Introduction}
Let $L$ be a Lagrangian submanifold in a \k manifold $(X,g,\omega)$. The second fundamental form $a_{ijk}$ is a symmetric 3-form on $L$. Let $H$ denote the mean curvature vector of $L$. Then the mean curvature 1-form $\alpha = \imath(H)\omega = a_ie^i$, where $a_i = a_{ijk}g^{jk}$, is a 1-form on $L$.\\
\begin{de}
A Lagrangian submanifold $L \subset (X,\omega)$ is called H-minimal (Hamiltonian minimal) if the mean curvature form $\alpha$ satisfies $d^*\alpha = 0$. $L$ is called L-minimal (Lagrangian minimal) if $\alpha$ is co-exact.\\
\end{de}
\begin{prop}
\label{aa}
An H-minimal (L-minimal) Lagrangian submanifold $L \subset (X,\omega)$ is a critical point of the volume functional restricted to the space of Hamiltonian (Lagrangian) deformations of $L$.
\end{prop}
{\bf Proof:} Let $F_t: L \rightarrow X$ be a Lagrangian (Hamiltonian) deformation family of Lagrangian submanifolds, $F_0 = {\rm id}$ and $V_t = \frac{dF_t}{dt}$. Then it is straightforward to derive that

\[
\left.\frac{d}{dt}\right|_{t=0} {\rm Vol}_{F_t^*g}(L) = -\int_L g(H,V_0)dV_{g|_L} = -\int_L g|_L(\alpha,\beta)dV_{g|_L},
\]

where $\beta = (\imath(V_0)\omega)|_L$. According to this first variation formula, the desired results are direct consequence of the fact that the deformation $F_t$ being Lagrangian (Hamiltonian) implies that $\beta$ is closed (exact). 
\hfill\rule{2.1mm}{2.1mm}\\

H-minimal Lagrangian submanifold in \k manifold was first investigated by Oh in \cite{O}. The H-minimal part of proposition \ref{aa} is essentially theorem 2.4 in \cite{O}, which was first observed by Weinstein according to Oh. The proof here basicly follow the argument in \cite{O}. In a recent work \cite{SW}, Schoen and Wolfson proved some important existence results for L-minimal Lagrangian surface in \k surface.\\

When the \k metric $g$ is K\"{a}hler-Einstein, the mean curvature form $\alpha$ is closed (\cite{D,B}). Namely, an H-minimal Lagrangian submanifold $L$ in a K\"{a}hler-Einstein manifold $(X,\omega)$ has harmonic mean curvature form (\cite{O}). If $\alpha$ is also exact, then $\alpha=0$ and $L$ is a minimal Lagrangian in $(X,\omega)$. In particular, if $(X,\omega)$ is Calabi-Yau, then the mean curvature form of any Lagrangian submanifold in $(X,\omega)$ is exact. Consequently, a Lagrangian submanifold in a Calabi-Yau manifold $(X,\omega)$ is H-minimal if and only if it is special. Therefore, H-minimal Lagrangian, which makes sense for any \k manifold, is both a generalization of special Lagrangian in Calabi-Yau manifold (\cite{HL}) and minimal Lagrangian in non-Calabi-Yau K\"{a}hler-Einstein manifold (\cite{B}).\\

Another important class of examples of H-minimal Lagrangian submanifolds are the top dimensional real torus orbits in a toric variety under a toric metric. The toric metric restricts to flat metric on such real torus orbit and the mean curvature form is constant, therefore co-close and actually harmonic. (The special cases of $\mathbb{C}^n$ with the flat metric and $\mathbb{C}P^n$ with the Fubini-Study metric are discussed in \cite{O}).\\ 

Special Lagrangian in Calabi-Yau manifolds and minimal Lagrangian in K\"{a}hler-Einstein manifolds are generally difficult to construct. Besides the obvious difficulty with the minimal surface equation, an important reason is our lack of understanding of the structure of Calabi-Yau metrics or more generally K\"{a}hler-Einstein metrics. In this work, we will develop methods to construct H-minimal Lagrangian fibration (theorem \ref{dg}) and minimal Lagrangian submanifold (theorem \ref{ec}) in K\"{a}hler-Einstein manifold through deformation from the toric model metric. This is made possible by H-minimal Lagrangian submanifolds, which give us additional flexbility to deform through general \k metrics that are neither K\"{a}hler-Einstein nor toric.\\

In \cite{ke2}, we discussed the degeneration of K\"{a}hler-Einstein metrics $\{g^{\rm KE}_t\}$ associated with a family of algebraic manifolds $\{X_t\}$ that degenerate into the central singular fibre $X_0$ under the so-called simple toroidal canonical degeneration. We proved in theorem 1.1 of \cite{ke2} that as $t$ approaches $0$, the K\"{a}hler-Einstein manifolds $\{(X_t,g^{\rm KE}_t)\}$ converge to the complete K\"{a}hler-Einstein manifold $(X_0\setminus {\rm Sing}(X_0),g^{\rm KE}_0)$ in the sense of Cheeger-Gromov. $X_0$ can be stratified into union of smooth equi-singular components. Points in $X_0$ that form the zero dimensional strata of such stratification will be called {\bf maximal degeneracy points} of $X_0$. Let $O\in X_0$ be one of such maximal degeneracy point of $X_0$. The following theorem is the main application of our construction.\\

{\bf Theorem \ref{ef}} {\it Let $O$ be a maximal degeneracy point in $X_0$. Then there exists a smooth family of minimal Lagrangian torus $L_t\subset (X_t,\omega^{\rm KE}_t)$ for $t$ small that approaches $O$ when $t$ approaches $0$.}\\
 
Deformation of H-minimal Lagrangian submanifolds in \k manifolds and its similarity to deformation of special Lagrangian submanifolds in \cy manifolds are discussed in section 2. In section 3, we construct H-minimal Lagrangian torus fibration for bounded perturbation of certain toric \k manifold through deformation method. In section 4, applying results from sections 2 and 3, we construct the minimal Lagrangian torus vanishing cycles in the toroidal degeneration family of \ke manifolds with negative first Chern class discussed in \cite{ke2}. Idea from \cite{sl1} is used in the construction to avoid singular deformation. (We notice the interesting construction of minimal Lagrangian tori in toric \ke manifolds with positive first Chern class by E. Goldstein \cite{G}, which may be viewed, in certain sense, as a dual situation.)\\

\se{Deformation of H-minimal Lagrangian in \k manifold}
Let $L$ be a Lagrangian submanifold in a \k manifold $(X,g,\omega)$. For $P\in L$, near $P$ locally we may choose holomorphic coordinate $z = x+iy$ such that $L=\{y=0\}$. Then $g_{i\bar{j}}$ are real along $L$. By adjusting $z$ by degree 2 polynomial on $z$, we may assume that $x$ is normal coordinate of $L$ at $P$ with respect to the Riemannian metric $g|_L$ on $L$. We will call such coordinate $z$ normal coordinate at $P$ for Lagrangian submanifold $L$ in \k manifold $X$. The computations in this section are all carried out under certain normal coordinate.\\

Let $\displaystyle e_i = \frac{\partial}{\partial x_i}$, $\displaystyle \tilde{e}_i = e_{\tilde{i}} = e_{n+i} = Je_i = \frac{\partial}{\partial y_i}$, and $\{e^i=dx_i,\tilde{e}^i = e^{\tilde{i}} = e^{n+i} = -Je^i = dy_i\}$ be the dual basis. Then the second fundamental form can be expressed as $a_{ijk} = (\nabla_i e_j, Je_k)$. It is straightforward to verify that at $P$, $\nabla_{e_i}e_j = a_{ij}^k \tilde{e}_k$, $\nabla_{e_i}\tilde{e}_j = -a_{ij}^k e_k$, $\nabla_{\tilde{e}_i}e_j = -a_{ij}^k e_k$.\\

Consider a family $\{L_t\}$ of Lagrangian submanifolds in $(X,\omega)$. We use the normal deformation vector fields $\{V_t\}$ for the family such that $V_t$ is orthogonal to $L_t$. Let $\beta_t = \imath(V_t)\omega$. Then we have\\
\begin{prop}
\label{ba}
\[
\dot{\alpha} = dd^*\beta - \imath(V_\beta){\rm Ric}^X,
\]

where $V_\beta$ is the vector field satisfying $g(V_\beta, W) = \beta(W)$.
\end{prop}
{\bf Proof:} It is straightforward to derive that
\[
\dot{g}_{ij} = (J\beta)_{i,j} + (J\beta)_{j,i}.
\]
Recall that
\[
\dot{\Gamma}_{ij}^k = \frac{1}{2}g^{kl}(\dot{g}_{il,j} + \dot{g}_{jl,i} - \dot{g}_{ij,l}).
\]
Hence
\[
\dot{\Gamma}_{ijk} = \dot{\Gamma}_{ij}^lg_{kl} = (J\beta)_{k,ij} + R_{jki}^l(J\beta)_l.
\]

Here the notation of the curvature is fixed by

\[
\beta_{k,ij} - \beta_{k,ji} = \beta_l R_{kij}^l = - \beta_l(R_{ijk}^l + R_{jki}^l).
\]

\[
\dot{a}_{ijk} = (\dot{\nabla}_i e_j, Je_k) = -\beta_{k,ij} + R_{j\tilde{l}\tilde{k}i}\beta^l.
\]

\[
\beta^X_{k,ij} = \beta^L_{k,ij} - (a_{ijl}a^{ls}_k + a_{jkl}a^{ls}_i)\beta_s.
\]
\[
R^X_{j\tilde{l}\tilde{k}i} = -R^X_{jlki} - 4{\rm Re}(R^X_{j\bar{l}k\bar{i}}).
\]

(There is a sign difference between the Riemannian and the \k curvature tensors.) Notice that at $P$, $g^X_{i\bar{j}} = g^L_{ij}$. Consequently $g_X^{i\bar{j}} = g_L^{ij}$. Hence

\[
R^X_{j\tilde{l}\tilde{k}i}g_L^{ik} = -4{\rm Re}(R^X_{j\bar{l}k\bar{i}})g_L^{ik} = -4{\rm Re}(R^X_{j\bar{l}k\bar{i}}g_X^{k\bar{i}}) = -4{\rm Re}(R^X_{j\bar{l}}) = -R^X_{jl}.
\]
\[
\beta^X_{k,ij}g_L^{ik} = \beta^L_{k,ij}g_L^{ik} - 2a_{jkl}a^{kls}\beta_s.
\]

When restricted to $L$,

\[
\dot{g}_{ij} = 2a_{ijk}\beta^k,\ \ \dot{g}^{ij} = -2a^{ijk}\beta_k.
\]
\[
\dot{a}_j = \dot{a}_{ijk}g^{ik} + a_{ijk}\dot{g}^{ik} = -\beta_{k,ij}g^{ik} - R^X_{jl}\beta^l = (d^*\beta)_j - R^X_{jl}\beta^l.
\]
\hfill\rule{2.1mm}{2.1mm}\\
\begin{prop}
\label{bb}
\[
\frac{d}{dt}(d^*\alpha) = D_\alpha\beta = d^*dd^*\beta - d^*(\imath(V_\beta){\rm Ric}^X) - V_\alpha(g(\alpha,\beta)) + 2(a^{ijk}a_j\beta_k)_i.
\]
\end{prop}
{\bf Proof:} 
\[
d^*\alpha = -{\rm Tr}(\nabla V_\alpha).
\]
\[
\frac{d}{dt}(d^*\alpha) = -{\rm Tr}(\dot{\nabla} V_\alpha) - {\rm Tr}(\nabla V_{\dot{\alpha}}) - {\rm Tr}(\nabla (\dot{g}^{ij}a_je_i)).
\]
\[
\dot{g}^{ij}a_j = -2a^{ijk}a_j\beta_k.
\]
\[
{\rm Tr}(\dot{\nabla} V_\alpha) = a^i\Gamma_{ij}^j = \frac{1}{2}a^ig^{jk}\dot{g}_{jk,i} = V_\alpha(g(\alpha,\beta)).
\]

Combining these calculations, we get the desired formula.
\hfill\rule{2.1mm}{2.1mm}\\

{\bf Remark:} Proposition \ref{bb} can also be derived from Oh's second variation formula (theorem 3.4 in \cite{O}).\\ 

Clearly, ${\rm Ker}(D_\alpha)$ represents the tangent space of the local deformation space of H-minimal Lagrangian submanifolds. It is straightforward to see that $D_\alpha: \Omega_{\rm closed}^1(L) \rightarrow \Omega^0_0(L)$ is an elliptic operator whose index equals to $h^1(L)$ (the first Betti number of $L$). Generically, when $D_\alpha$ is surjective, the deformation space of H-minimal Lagrangian near $L$ is smooth of dimension $h^1(L)$, which is a complete analog to the case of special Lagrangian (\cite{ML}). On the other hand, unlike the case of special Lagrangian, in general, the deformation of $L$ could be obstructed and local deformation space could be of higher dimension than $h^1(L)$ in the case of H-minimal Lagrangian. A good example is the case of Riemann sphere $S^2$ with the standard round metric. H-minimal Lagrangians are exactly round circles of constant mean curvature in $S^2$ and minimal Lagrangians are the great circles in $S^2$. The dimensions of the moduli spaces of both are greater than $1$. The reason is that the round metric on $S^2$ is a very special metric. A choice of more generic metric on $S^2$ like that of American football will result in 1-parameter family of H-minimal Lagrangian circles genericly. In this paper, we will not need to consider the non-generic situations.\\

\se{H-minimal Lagrangian fibration}
In this section, we will construct H-minimal Lagrangian torus fibration for bounded perturbation of certain toric \k manifold (considered in \cite{ke2,kesl}) through deformation method. One key idea that makes the deformation possible is that the bounded perturbation of the toric \k metric we consider can be reduced to a small perturbation of another toric \k metric (proposition \ref{db}). We will start with the formulism in \cite{kesl}. Let $F: (\mathbb{C}^*)^n \rightarrow \mathbb{R}^n$ be defined as $x = F(z) = (\log |z_1|^2,\cdots,\log |z_n|^2)$. Then for any bounded convex set $\Delta \subset \mathbb{R}^n$, one can define the generalized cylinder $D_\Delta = F^{-1}(\Delta)$. \\

A convex polyhedron can be defined through an equivalence class (modulo linear functions) of convex piecewise linear integral functions $w$ on a lattice $M$ that are compatible with a complete fan $\Sigma$ in $M$. Let $\Sigma(k)$ denote the set of $k$-dimensional cones in $\Sigma$. Assume $\Sigma$ is rational and simplicial. Then $\Sigma(1)$ can also be identified as a subset of $M$ containing the primitive integral elements of the corresponding 1-dimensional cones in $\Sigma(1)$. We may write $w = \{w_m\}_{m\in \Sigma(1)}$, where $w_m$ is the value of $w$ at $m\in \Sigma(1)$. From these data, one may define convex polyhedrons

\[
\Delta_\tau = \tau \Delta,\ \Delta = \{x \in N_{\mathbb{R}}|\langle m,x\rangle + w_m \leq 0,\ {\rm for}\ m\in \Sigma(1)\}.
\]

Let

\[
\rho_\tau(x) = \rho(x/\tau) - n\log \tau^2,\ \rho(x) = - \sum_{m\in \Sigma(1)} \log \left(\langle m,x\rangle + w_m\right)^2.
\]

\[
\omega_\tau = \partial \bar{\partial} \rho_\tau(x) = \sum_{m\in \Sigma(1)} \frac{\partial q_m \bar{\partial} q_m}{q_m^2},\ \ {\rm where}\ q_m = w_m + \frac{1}{\tau}\langle m,x\rangle
\]

defines a complete toric \k metric on the generalized cylinder $D_{\Delta_\tau}$. $\rho$ is strictly convex. We fix the origin to be the unique critical point of $\rho$. Then clearly the origin is also the unique critical point of $\rho_\tau$ for all $\tau$. The only properties of $\rho_\tau$ we will need are the following\\
\begin{equation}
\label{ca}
\rho_\tau(x) = \rho(x/\tau) + C(\tau),\ \ \lim_{c\rightarrow 1} \rho|_{\Delta_c} = +\infty.
\end{equation}

Let $\hat{g}_\tau = \tau^2g_\tau$, we have\\
\begin{lm}
\label{cb}
For any $c\in (0,1)$, $\hat{g}_\tau|_{F^{-1}(x)}$ is a flat metric with bounded geometry for $x \in \Delta_{c\tau}$. (The bound depends on $c$ and is uniform on $x \in \Delta_{c\tau}$.)
\end{lm}
{\bf Proof:} (\ref{cb}) implies that $\hat{g}_\tau|_{F^{-1}(x)} = \hat{g}_1|_{F^{-1}(x/\tau)}$. Namely, the lemma can be reduced to the special case of $\tau=1$, which is quite obvious.
\hfill\rule{2.1mm}{2.1mm}\\
\begin{lm}
\label{da}
Let $(T,h)$ be a torus with flat metric $h$ of bounded geometry. Assume that a function $f$ on $T$ is $C^\infty$-bounded with respect to $h_\tau = \tau^{-2}h$ and satisfies $\displaystyle \int_TfdV_h =0$. Then for each positive integer $n$, there exists a constant $C(n)>0$ independent of $\tau$ such that $|f| \leq C(n)\tau^{-n}$.
\end{lm}
{\bf Proof:} Since $\displaystyle \int_TfdV_h =0$, there exists $\theta_0\in T$ such that $f(\theta_0)=0$. Hence

\[
|f(\theta)| = |f(\theta)-f(\theta_0)| \leq C|\nabla f|_{h_\tau}{\rm Diam}(T,h_\tau) \leq C\tau^{-1}.
\]

This estimate can also be rewritten as

\[
|f| \leq C\left|\tau\frac{\partial f}{\partial \theta}\right|\tau^{-1} \leq C\tau^{-1}.
\]

It is easy to observe that $\tau\frac{\partial f}{\partial \theta}$ will satisfy all the assumptions for $f$ in the lemma. By induction, we get the conclusion of the lemma.
\hfill\rule{2.1mm}{2.1mm}\\

(Two quasi-isometric \k forms $\omega$ and $\omega' = \omega + \partial \bar{\partial} f$ are called $C^\infty$-quasi-isometric if $f$ is $C^\infty$-bounded with respect to $\omega$. Consequently, $\omega'$ is a $C^\infty$-bounded tensor with respect to $\omega$.) Consider a family of \k metrics $\omega'_\tau = \omega_\tau + \partial \bar{\partial} f_\tau$. Let $\hat{\omega}'_\tau = \tau^2\omega'_\tau$.\\
\begin{prop}
\label{db}
Assume that $\omega'_\tau$ is $C^\infty$-quasi-isometric to $\omega_\tau$ (uniform with respect to $\tau$). There exists a decomposition $\hat{\omega}'_{\tau} = \hat{\omega}^0_{\tau} + \hat{\omega}^1_{\tau}$ such that $\hat{\omega}^0_{\tau} = \hat{\omega}_\tau + \tau^2\partial \bar{\partial} f^0_\tau$ is toric and $\hat{\omega}^1_{\tau} = \tau^2\partial \bar{\partial} f^1_\tau$. For any $c\in (0,1)$, index set $I$ and positive integer $n$, there exists a constant $C(n,c,I)>0$ independent of $\tau$ such that $|\nabla_I f^1_\tau| \leq C(n,c,I)\tau^{-n}$ in $D_{\Delta_{c\tau}}$ with respect to $\hat{\omega}_{\tau}$.
\end{prop}
{\bf Proof:} There is a canonical decomposition $f_\tau = f^0_\tau + f^1_\tau$ such that $f^0_\tau$ is constant in each fibre $F^{-1}(x)$ and the integral of $f^1_\tau$ on each fibre $F^{-1}(x)$ is zero. In another word, $f^0_\tau$ is the average function of $f_\tau$ along fibres of $F$. This gives us the desired decomposition $\hat{\omega}'_{\tau} = \hat{\omega}^0_{\tau} + \hat{\omega}^1_{\tau}$.\\

Lemma \ref{cb} implies that for any $c\in (0,1)$, $F^{-1}(x)$ is of bounded geometry with respect to $\hat{\omega}_\tau$ for $x\in \Delta_{c\tau}$. (The bound depends on $c$.) Apply lemma \ref{da} to $\nabla_I f^1_\tau$, we get the desired estimate.
\hfill\rule{2.1mm}{2.1mm}\\
\begin{co}
\label{dd}
For any $c\in (0,1)$, index set $I$ and positive integer $n$, there exists a constant $C(n,c,I)>0$ independent of $\tau$ such that $|\nabla_I f^1_\tau| \leq C(n,c,I)\tau^{-n}$ in $D_{\Delta_{c\tau}}$ with respect to ${\omega}_{\tau}$. Consequently, if ${\omega}'_{\tau}$ is K\"{a}hler-Einstein, then ${\omega}^0_{\tau}$ is K\"{a}hler-Einstein up to $O(1/\tau)$-perturbation.
\hfill\rule{2.1mm}{2.1mm}\\
\end{co}
Define a family of metric $\hat{\omega}_{\tau,s} = \hat{\omega}^0_\tau + s\tau^2\partial \bar{\partial} f^1_\tau$. Then $\hat{\omega}_{\tau,0} = \hat{\omega}^1_\tau$ and $\hat{\omega}_{\tau,1} = \hat{\omega}'_\tau$. Let $V$ denote the Hamiltonian-gradient vector field (see \cite{sl1} and references therein) associated with the family $\{\hat{\omega}_{\tau,s}\}_{s\in [0,1]}$ of \k metrics, and $\phi_s$ be the corresponding Hamiltonian-gradient flow. Then $\phi_s^*\hat{\omega}_{\tau,s} = \hat{\omega}_{\tau,0} = \hat{\omega}^0_\tau$.\\ 
\begin{lm}
\label{dc}
$V = -\tau^2\nabla_{\tau,s}f^1_\tau$ (the gradient vector field of $-\tau^2f^1_\tau$ with respect the \k metric $\hat{g}_{\tau,s}$). For any $c\in (0,1)$, index set $I$ and positive integer $n$, there exist positive constants $C(n,c,I),C'(n,c,I)$ independent of $\tau$ such that $|\nabla_I V| \leq C(n,c,I)\tau^{-n}$, $|\nabla_I (\phi^*_s\hat{g}_{\tau,s} - \hat{g}_{\tau,0})| \leq C'(n,c,I)\tau^{-n}$ in $D_{\Delta_{c\tau}}$ with respect to $\hat{\omega}_{\tau}$.
\end{lm}
{\bf Proof:} $2{\rm Re}\left(\frac{\partial}{\partial s}\right) + V$ is perpendicular to any vector field $W$ on $D_{\Delta_{c\tau}}$ with respect to the following \k metric on $(z,s)$

\[
\hat{\omega}^0_\tau + \tau^2\partial\bar{\partial} ({\rm Re}(s)f^1_\tau) = \hat{\omega}_{\tau,s} + \frac{\tau^2}{2}(ds\bar{\partial} f^1_\tau + \partial f^1_\tau d\bar{s}).
\]

Hence

\[
\hat{g}_{\tau,s}(V,W) = -\tau^2\langle df^1_\tau,W\rangle,\ \ V = -\tau^2\nabla_{\tau,s}f^1_\tau.
\]

With this expression of $V$, the rest of the lemma is a consequence of corollary \ref{dd}.
\hfill\rule{2.1mm}{2.1mm}\\

Let $L_0 = F^{-1}(x)$. Then $\phi_s(L_0)$ is a Lagrangian torus with respect to $\hat{\omega}_{\tau,s}$. Let $(\theta,y)$ be the toric Darboux coordinate with respect to the toric metric $\hat{\omega}^0_\tau$ such that $y|_{L_0} = 0$. For a function $h(\theta)$ on $L_0$, let $L(h)$ be the graph of $y=dh(\theta)$ in the symplectic neighborhood of $L_0$. $\Phi(h,s) = d^*\alpha_{L(h)}$ defines a map $\Phi: {\cal B}_1\times\mathbb{R} \rightarrow {\cal B}_2$, where ${\cal B}_1 = C^{4,\alpha}_0(L_0)$, ${\cal B}_2 = C^{\alpha}_0(L_0)$ and $\alpha_{L(h)}$ is the mean curvature form of $L(h)$ under the \k metric $\hat{\omega}_{\tau,s}$.\\
\begin{lm}
\label{de}
\[
\left\|\frac{\partial \Phi}{\partial h}(h,s) - \frac{\partial \Phi}{\partial h}(0,0)\right\| = O(\tau^{-1},|h|_{{\cal B}_1}),\ \ \|\Phi(0,s)\|_{{\cal B}_2} = O(\tau^{-1}).
\]
\end{lm}
{\bf Proof:} The estimates

\[
\left\|\frac{\partial \Phi}{\partial h}(h,s) - \frac{\partial \Phi}{\partial h}(h,0)\right\| = O(\tau^{-1}),\ \ \|\Phi(0,s)\|_{{\cal B}_2} = O(\tau^{-1})
\]

are easy consequences of lemma \ref{dc}. The estimate

\[
\left\|\frac{\partial \Phi}{\partial h}(h,0) - \frac{\partial \Phi}{\partial h}(0,0)\right\| = O(|h|_{{\cal B}_1})
\]

is straightforward to derive.
\hfill\rule{2.1mm}{2.1mm}\\
\begin{lm}
\label{df}
\[
\left\|\left(\frac{\partial \Phi}{\partial h}\right)^{-1}(0,0)\right\| \leq C.
\]
\end{lm}
{\bf Proof:} 
\[
\frac{\partial \Phi}{\partial h}(0,0)\delta h = D_\alpha d\delta h.
\]

With respect to the rescaled metric $|{\rm Ric}^X| = O(1/\tau^2)$, $|a| = O(1/\tau)$ and $|\alpha| = O(1/\tau)$. Therefore

\[
\frac{\partial \Phi}{\partial h}(0,0)\delta h = \Delta^2\delta h + O(1/\tau^2).
\]

Since $\|\Delta^{-2}\|\leq C$, for $\tau$ large, we have the desired estimate.
\hfill\rule{2.1mm}{2.1mm}\\
\begin{theorem}
\label{dg}
Assume the \k potential $\rho_\tau$ of the toric \k metric $\omega_\tau$ satisfies (\ref{ca}) and the \k metric $\omega'_\tau = \omega_\tau + \partial \bar{\partial} f_\tau$ is $C^\infty$-quasi-isometric to $\omega_\tau$ (uniform with respect to $\tau$). Fix $c\in (0,1)$, when $\tau$ is large enough, there exist a smooth family of H-minimal Lagrangian torus fibration $F_s$ over $D_{\Delta_{c\tau}}$ with respect to the \k form and metric $(\hat{\omega}^0_\tau, \phi^*_s\hat{g}_{\tau,s})$. $F_s$ are $\tau^{-1}$-perturbations of the toric fibration $F_0 = F$. When $x\in \Delta_{c\tau}$ varies, $\phi_1(F_1^{-1}(x))$ forms an H-minimal Lagrangian torus fibration under the \k form $\omega'_{\tau}$.
\end{theorem}
{\bf Proof:} Lemmas \ref{de} and \ref{df} enable us to apply the quantitative implicit function theorem (theorem 3.2 in \cite{sl1}) to $\Phi$. Consequently, there exist a constant $C_1$ and a unique family $\{h_s\}_{s\in [0,1]}$ such that $h_0=0$, $|h_s|_{{\cal B}_1} \leq C_1$ and $L(h_s)$ is an H-minimal Lagrangian with respect to the \k form and metric $(\hat{\omega}^0_\tau, \phi^*_s\hat{g}_{\tau,s})$. Further more, $h_s$ actually satisfies $|h_s|_{{\cal B}_1} = O(1/\tau)$.\\

To show that $L(h_s)$ forms a fibration when $L_0$ varies, notice that from the estimate $|h_s|_{{\cal B}_1} = O(1/\tau)$, we have that $L(h_s)$ is an $O(1/\tau)$-perturbation of $L_0$. Since $L_0$ is toric, the metric on $L_0$ is flat and the second fundamental form on $L_0$ is constant. According to proposition \ref{bb}, the H-minimal Lagrangian deformation 1-forms are exactly the constant 1-forms on $L_0$. Consequently, the H-minimal Lagrangian deformation 1-forms on $L(h_s)$ are $O(1/\tau)$-perturbations of the constant 1-forms, therefore are non-vanishing anywhere. This implies that $L(h_s)$ forms a fibration when $L_0$ varies. We take this fibration to be $F_s$.
\hfill\rule{2.1mm}{2.1mm}\\

\se{Minimal Lagrangian torus in K\"{a}hler-Einstein manifold}
\begin{lm}
\label{eb}
The logrithm of the volume of $F^{-1}(x)$ under a toric \k metric $\omega$ forms a function $u(x)$. If the Ricci curvature of $\omega$ is negative, then the $x$ where $u(x)$ reaches minimal if exists will be unique.
\end{lm}
{\bf Proof:} The condition of the lemma implies that $u$ is a strictly convex function of $x$. Therefore, the minimal if exists will be unique.
\hfill\rule{2.1mm}{2.1mm}\\

We will call $u(x)$ in lemma \ref{eb} the logrithm of the volume function under the toric metric $\omega$. Recall that ${\omega}'_{\tau} = {\omega}^0_{\tau} + {\omega}^1_{\tau}$ such that ${\omega}^0_{\tau} = {\omega}_\tau + \partial \bar{\partial} f^0_\tau$ is toric. Let $u^0_\tau$ denote the logrithm of the volume function under ${\omega}^0_{\tau}$. We have\\
\begin{lm}
\label{ea}
There exist $c\in (0,1)$ and a unique $x_0 \in \Delta_{c\tau}$ so that $u^0_\tau$ reaches the minimal at $x_0$ and $|u^0_\tau(x) - u^0_\tau(x_0)| \geq 1$ for $x \in \partial \Delta_{c\tau}$.
\end{lm}
{\bf Proof:} (\ref{ca}) implies that for any $C>0$ there exist $c\in (0,1)$ such that $|\rho_\tau(x) - \rho_\tau(0)| = |\rho(x/\tau) - \rho(0)| \geq C$ for $x \in \partial \Delta_{c\tau}$.\\

The \k potential of ${\omega}^0_{\tau}$ is $\rho^0_\tau = \rho_\tau + f^0_\tau$. Since $f^0_\tau$ is uniformly bounded, there exists $c\in (0,1)$ such that $|\rho^0_\tau(x) - \rho^0_\tau(0)| \geq 2$ for $x \in \partial \Delta_{c\tau}$. Since ${\omega}^0_{\tau}$ is \ke up to $O(1/\tau)$-perturbation according to corollary \ref{dd}, ${\omega}^0_{\tau}$ clearly has negative Ricci curvature, also $\rho^0_\tau$ is a $O(1/\tau)$-perturbation of the logrithm of the volume function $u^0_\tau$. When $\tau$ is large, $|u^0_\tau(x) - u^0_\tau(0)| \geq 1$ for $x \in \partial \Delta_{c\tau}$. Lemma \ref{eb} implies that there exist unique $x_0 \in \Delta_{c\tau}$, where $u^0_\tau$ reaches the minimal and $|u^0_\tau(x) - u^0_\tau(x_0)| \geq 1$ for $x \in \partial \Delta_{c\tau}$.
\hfill\rule{2.1mm}{2.1mm}\\

The moduli space of Lagrangian torus in a symplectic manifold modulo Hamiltonian equivalence locally around a Lagrangian torus $L_0$ can be naturally identified with $H^1(L_0)$. A Lagrangian torus fibration in a \k manifold is called {\bf closed} if the mean curvature form of each Lagrangian torus fibre is a closed 1-form. A closed Lagrangian torus fibration near $L_0$ naturally induces a map $\Psi: H^1(L_0) \rightarrow H^1(L_0)$ defined as $\Psi([L]) = [\alpha_L]$. Proposition \ref{ba} implies that\\
\begin{prop}
\label{ed}
The tangent map $d\Psi: H^1(L) \rightarrow H^1(L)$ has the expression
\[
d\Psi([\beta]) = -[\imath(V_\beta){\rm Ric}^X] 
\]
\hfill\rule{2.1mm}{2.1mm}
\end{prop}
A Lagrangian torus fibration in a K\"{a}hler-Einstein manifold is automatically closed. Another important example of closed Lagrangian torus fibration is the toric torus fibration under a toric \k metric. Since ${\omega}^0_{\tau}$ is toric, if ${\omega}'_{\tau}$ is K\"{a}hler-Einstein, the H-minimal Lagrangian fibrations $F_0=F$ and $F_1$ in theorem \ref{dg} are both closed Lagrangian fibrations. They induce maps $\Psi_0$ and $\Psi_1$.\\ 
\begin{theorem}
\label{ec}
Assume the \ke metric $\omega'_\tau = \omega_\tau + \partial \bar{\partial} f_\tau$ is $C^\infty$-quasi-isometric to $\omega_\tau$ (uniform with respect to $\tau$) and the \k potential $\rho_\tau$ of the toric \k metric $\omega_\tau$ satisfies (\ref{ca}). For suitable $c\in (0,1)$, when $\tau$ is large enough, there exists a unique $x_1\in \Delta_{c\tau}$ so that $\phi_1(F_1^{-1}(x_1))$ is a minimal Lagrangian torus under the \k form $\omega'_{\tau}$.
\end{theorem}
{\bf Proof:} Proposition \ref{ed} and corollary \ref{dd} imply that $d\Psi_1 = {\rm id}$, $d\Psi_0 = {\rm id} + O(1/\tau)$ and $|\Psi_1(L) - \Psi_0(L)|_{\hat{g}^0_{\tau}} = O(1/\tau)$, where $L = F_0^{-1}(x)$ for some $x\in \Delta_{c\tau}$. (It is conceptually more clear to use the rescaled metric $\hat{g}^0_{\tau}$, under which $L$ is of bounded geometry.) Since $\Psi_0(L_0)=0$ for the fibre $L_0 = F_0^{-1}(x_0)$ with minimal volume under ${\omega}^0_{\tau}$, there exist a unique $L_1 = F_0^{-1}(x_1)$ such that $\Psi_1([L_1]) = 0$ and $|[L_1] - [L_0]|_{\hat{g}^0_{\tau}} = O(1/\tau)$. Consequently, the mean curvature 1-form of $F_1^{-1}(x_1)$ is closed, coclosed and exact, therefore vanishes. When $\tau$ is large, $|u^0_\tau(x) - u^0_\tau(x_1)| \geq 1 - O(1/\tau) >0$ for $x \in \partial \Delta_{c\tau}$. Namely $x_1 \in \Delta_{c\tau}$. 
\hfill\rule{2.1mm}{2.1mm}\\

We are now ready to discuss the main application of our results as mentioned in the introduction. We will follow the notations as in the introduction. Using the Hamiltonian-gradient flow with respect to the parameter $t$, it is straightforward to show that the vanishing cycle in $X_t$ that vanishes to $O\in X_0$ when $t$ approaches $0$ can be represented by Lagrangian torus in $X_t$. We will show in the following theorem that such vanishing cycle can actually be represented by a minimal Lagrangian torus in $(X_t,g^{\rm KE}_t)$. As in the construction of Lagrangian representative, one may consider constructing the minimal Lagrangian representitive of the vanishing cycles through deformation starting from $O$. Such method will run into singular deformation problem that is usually very delicate to handle if it is solvable at all. Instead, in our construction we will use an idea similar to the key idea in \cite{sl1}, by constructing suitable local model (family of model \k manifolds), where the solution is clear, then deform to the actual K\"{a}hler-Einstein manifold $(X_t,g^{\rm KE}_t)$ for fixed $t(\not=0)$ small, therefore avoiding the singular deformation problem.\\

Around a maximal degeneracy point $O\in X_0$, the total space ${\cal X}$ is locally toric. $X_t$ locally is the image of the toric embedding $\{s_m = e^{\tau w_m/2}z^m\}_{m\in \Sigma(1)}$, where $\tau = -\log |t|^2$. Under coordinate $z$, $X_t$ locally near $O$ can be identified with $F^{-1}(\Delta_{\tau})$.
\begin{prop}
\label{ee}
For certain fixed $\mu>0$, in $F^{-1}(\Delta_{\tau-\mu})$, the K\"{a}hler-Einstein metric on $X_t$ can be expressed as

\[
\omega^{\rm KE}_t = \partial \bar{\partial} \rho^{\rm KE}_\tau,\ \ \rho^{\rm KE}_\tau = b^0 + b^1 - \sum_{m\in \Sigma(1)} \log \left(\tau q_m + b_m\right)^2,
\]

where $b^0$ and $b_m$ for $m\in \Sigma(1)$ are $C^\infty$-bounded functions of $\{s_m\}_{m\in \Sigma(1)}$, and $b^1$ is a $C^\infty$-bounded function with respect to $\omega_\tau = \partial \bar{\partial} \rho_\tau$.
\end{prop}
{\bf Proof:} In this proof, we will use notations from \cite{ke2}. Definition of $\|s^I_m\|_m$ in section 2 of \cite{ke2} implies that for certain fixed $\mu>0$, in $F^{-1}(\Delta_{\tau-\mu})$, $\log \|s^I_m\|_m = \tau q_m + b_m$, where $b_m$ for $m\in \Sigma(1)$ are $C^\infty$-bounded functions of $\{s_m\}_{m\in \Sigma(1)}$. Formula (3.1) of \cite{ke2} implies that the \k potential $b^0$ of $\hat{\omega}_t$ can be made a $C^\infty$-bounded function of $\{s_m\}_{m\in \Sigma(1)}$. Definition of the approximate metric $\omega_t$ in section 4 of \cite{ke2} implies that $\displaystyle b^0 - \sum_{m\in \Sigma(1)} \log \left(\tau q_m + b_m\right)^2$ is a \k potential of $\omega_t$. Assume the \ke metric $\omega^{\rm KE}_t = \omega_t + \partial \bar{\partial} b^1$. The Monge-Amp\`{e}re estimate of the \ke metric (\cite{yau}) together with the estimate in proposition 4.6 of \cite{ke2} imply that $b^1$ is a $C^\infty$-bounded function with respect to $\omega_\tau = \partial \bar{\partial} \rho_\tau$.
\hfill\rule{2.1mm}{2.1mm}\\
\begin{theorem}
\label{ef}
Let $O$ be a maximal degeneracy point in $X_0$. Then there exists a smooth family of minimal Lagrangian torus $L_t\subset (X_t,\omega^{\rm KE}_t)$ for $t$ small that approaches $O$ when $t$ approaches $0$.
\end{theorem}
{\bf Proof:} On $X_t$, apply proposition \ref{ee}, we have

\[
\rho^{\rm KE}_\tau = \rho_\tau + f_\tau,\ \ f_\tau = b^0 + b^1 - \sum_{m\in \Sigma(1)} \log \left(1 + \frac{b_m}{\tau q_m}\right)^2.
\]

Since $\mu$ is fixed, for a fixed $c\in (0,1)$, when $\tau$ is large enough, we have $\Delta_{c\tau} \subset \Delta_{\tau - \mu}$. In $F^{-1}(\Delta_{c\tau})$, $|s_m|^2 \leq e^{-(1-c)\tau}$ and $|\log |s_m|^2| \geq (1-c)\tau$. A proper basis (as defined in the proof of proposition 4.5 of \cite{ke2}) for $\omega_\tau$ in $F^{-1}(\Delta_{c\tau})$ is $\displaystyle \left\{\tau z_i\frac{\partial}{\partial z_i}\right\}_{i=1}^n$. Since

\[
\tau z_i\frac{\partial s_m}{\partial z_i}= \tau m_is_m = O(\tau e^{-(1-c)\tau/2}),
\]
\[
\frac{1}{q_m}= O\left(1\right),\ \ \tau z_i\frac{\partial}{\partial z_i}\left(\frac{1}{q_m}\right) = -\frac{m_i}{q_m^2}= O\left(1\right).
\]

When $\tau$ is large, $s_m$ and $\displaystyle \frac{1}{q_m}$ for $m\in \Sigma(1)$ are $C^\infty$-bounded functions in $F^{-1}(\Delta_{c\tau})$ with respect to $\omega_\tau$. Therefore, $b^0$ and $b_m$ for $m\in \Sigma(1)$, which are $C^\infty$-bounded functions of $\{s_m\}_{m\in \Sigma(1)}$ according to proposition \ref{ee}, are $C^\infty$-bounded functions in $F^{-1}(\Delta_{c\tau})$ with respect to $\omega_\tau$. Consequently, $f_\tau$, which is a function of $b^0$, $b^1$, $\displaystyle \frac{1}{q_m}$ and $b_m$ for $m\in \Sigma(1)$, is a $C^\infty$-bounded function in $F^{-1}(\Delta_{c\tau})$ with respect to $\omega_\tau$. (Notice that the bounds of $f_\tau$ is independent of $c$ and $\tau$ as long as $\tau$ is taken to be suitably large according to $c$.) Apply theorem \ref{dg} and theorem \ref{ec}, we get the desired minimal Lagrangian torus $L_t\subset X_t$ for $t$ small.
\hfill\rule{2.1mm}{2.1mm}\\

\ifx\undefined\bysame
\newcommand{\bysame}{\leavevmode\hbox to3em{\hrulefill}\,}
\fi


\noindent
\begin{thebibliography}{100}
\bibitem{B}
Bryant, R.,
{\em Minimal Lagrangian submanifolds of \ke manifolds},
Lecture Notes in Mathematics, 1255, 1985.
%
\bibitem{D}
Dazord, P.,
{\em Sur la g\'{e}ometrie des sous-fibr\'{e}s et des feuilletages lagrangiense}, Ann. Sci. \'{E}c. Norm. Super., IV. S\'{e}r. 13, (1981), 465-480.
%
\bibitem{G}
Goldstein, E., {\em Minimal Lagrangian tori in \ke manifolds}, math.DG/0007135.
%
\bibitem{GT}
Gilbarg, D. and Trudinger, N.S.,
{\em Elliptic partial differential equations of second order},
Second edition, Springer-Verlag, 1983.
%
\bibitem{HL}
Harvey, R. and Lawson, H.B.,
{\em Calibrated Geometries},
Acta Math. 148 (1982), 47-157.
%
\bibitem{ML}
Mclean, R., {\em Deformations of calibrated submanifolds},
Comm. Anal. Geom. 6 (1998), 705-747.
% 
\bibitem{O1}
Oh, Y., {\em Second variation and stabilities of minimal Lagrangian submanifolds},
Invent. Math. 101 (1990), 501-519.
% 
\bibitem{O}
\bysame, {\em Volume minimization of Lagrangian submanifolds under Hamiltonian deformations},
Math. Z. 212 (1993), 175-192.
% 
\bibitem{ke2}
Ruan, W.-D., {\em Degeneration of K\"{a}hler-Einstein metrics II}, preprint, math.DG/0303113.
% 
%
\bibitem{sl1}
\bysame,
{\em Generalized special Lagrangian torus fibrations of Calabi-Yau hypersurfaces I}, preprint, math.DG/0303114.
%
% 
\bibitem{kesl}
\bysame,
{\em Degeneration of complete toric K\"{a}hler-Einstein manifolds and minimal Lagrangians},
in preparation.
%
\bibitem{SW}
Schoen, R. and Wolfson, J.,
{\em Minimizing area among Lagrangian surfaces: the mapping problem},
J. Differential Geom. 58 (2001), 1-86.
%
\bibitem{yau}
Yau, S.-T., {\em On the Ricci curvature of a compact K\"{a}hler manifold and the complex Monge-Amp\`{e}re equation, I}, Comm. Pure. and Appl. Math.,  {\bf 31} (1978), 339-411.
\end{thebibliography}
\end{document}